# A Note on Computation of Number of Fuzzy Bitopological Space

**Bhimraj Basumatary**
Assistant Professor, Department of Mathematical Sciences, Bodoland University, BTAD, India

**Jili Basumatary**
Research Scholar, Department of Mathematical Sciences, Bodoland University, BTAD, India

**Nijwm Wary**
Research Scholar, Department of Mathematical Sciences, Bodoland University, BTAD, India

**Abstract:** In this article we compute the number of fuzzy bitopological space with having two open sets, three open sets, four open sets and five open sets. Also, we have given some results on computation of number of fuzzy bitopological space.

1. **Introduction:**

A fuzzy set is defined by Lotfi A. Zadeh [1], is a function from a set $X$ to $[0,1]$. This definition has been extended to a more general set than the unit interval; for example to a complete lattice. General topology was one of the first branches of pure mathematics to which fuzzy sets have been applied systemically. Starting from single topology it is extended to bitopology and tritopology etc. V. Krishnamurthy [2] computed some results on the number of topologies on a finite set. Where he obtained a sharper bound for number of distinct topology, namely, $2^{2^n}$. H. Sharp [3] shown that no topology, other than the discrete, has cardinal greater than $\frac{3}{4}2^n$ and some other bounds are derived on the cardinality of connected, non-$T_0$, connected and non-$T_0$, and non-connected topologies. Richard P. Stanley[4] found all non-homeomorphic topologies with n-points and $>\frac{7}{16}2^n$ open set by using the correspondence between finite $T_0$-topologies and partial orders. He computed number of topology on a finite set having $n$ elements and $k$ open set i.e. $T(n,k)$ for large values of $k$, viz.; $3.2^{n-3} < k < 2^n$ and also $T_0$ topologies on a set having either $n+1, n+2$ or $n+3$ open sets. G.A. Kamel [5] formulated special case for computing the number of chain topological spaces, and maximal elements with the natural generalization. Kari Ragnarsson, Bridget Eileen Tenner [6], studied the smaller possible number of points in a topological space having $k$-open sets. Moussa Benoumhani [7] computed number of topology on a finite set $X$ having $n$ elements and $k$ open set for $2 \leq k \leq 12$, as well as other results concerning $T_0$ topologies on $X$ having $n+4 \leq k \leq n+6$ open sets.

In 1968, C.L. Chang [8], introduced the concept of fuzzy topological spaces as an application of fuzzy sets to general topological spaces. Many researchers has developed various properties of topology and they extended general topology to fuzzy topology. M. Benoumhani and M. Kolli [9] also studied about fuzzy topologies and partitions. A. Jaballah and F.B. Saidi [10] studied about the Length of maximal chains and number of ideals in commutative rings. M. Benoumhani and A. Jaballah [11] also studied about fuzzy topological spaces and computed some results for finding number of fuzzy topology.

Throughout this paper we use the following symbols:

Let $X$ be a set having $n$ elements, $M$ be a totally ordered one having $m$ elements and let $\mathcal{F}$ be the collection of all fuzzy subsets of $X$ with membership values in $M$.

$\mathcal{F}$ is partially ordered by $\mu \leq \vartheta \Leftrightarrow \mu(x) \leq \vartheta(x)$ for every $x \in X$. This set is a complete lattice with the same partial order. We also have $\mu < \vartheta$, if and only if $\mu \leq \vartheta$ and $\mu(x) < \vartheta(x)$, for some $x \in X$.

The fuzzy subsets $0_\mathcal{F}(x) = 0$ for every $x \in X$, and $1_\mathcal{F}(x) = 1$ for every $x \in X$. For every fuzzy subset $\mu$ different from $0_\mathcal{F}$ and $1_\mathcal{F}$, we have $0_\mathcal{F} < \mu < 1_\mathcal{F}$.

In this paper we compute the number of fuzzy bitopological spaces having same number of open sets $k$, where $k$ is small, in each pair of fuzzy topology on a finite set $X$ having $n$ element whose membership value lies in $M$ having $m$ elements.

By $\tau_\mathcal{F}(n, m, k)$ we denote number of fuzzy topology on a set $X$ having $n$ elements, here $m$ is the number of elements in $M$ and $k$ is the number of open set. Also by $(\tau_i, \tau_j)_\mathcal{F}(n, m, k)$ we denote number of fuzzy bitopological space where $m, n$ and $k$ are as defined above.

**2. Preliminary:**

**Fuzzy set:**

**Definition 2.1:** Let $X$ be a universal set. Then a function $A: X \to [0,1]$ define a fuzzy set on $X$ where $A$ is called the membership function and $\mu_A(x)$ is called membership grade of $x$.

We also write fuzzy set

$A = \{(x, \mu_A(x)): x \in X\}$, Where each pair $(x, \mu_A(x))$ is called singleton.

**Definition 2.2:** A fuzzy topology $\tau$ on a set $X$ consist of a collection of fuzzy subsets of $X$ called open sets, satisfying the following three axioms:

(1) The fuzzy subsets $0_\mathcal{F}$ and $1_\mathcal{F}$ are in $\tau$.
(2) The union $\vee_{i \in I} u_i$ of any collection $\{u_i : i \in I\}$ of elements of $\tau$ is also in $\tau$.
(3) The intersection $u_1 \cap u_2$ of any two elements $u_1$ and $u_2$ of $\tau$ is also in $\tau$.

The existence of a topology $\tau$ in the collection of fuzzy subsets of $\mathcal{F}$ with membership values in $M$ implies necessarily that $0_\mathcal{F}$ and $1_\mathcal{F}$ are open sets in $\tau$.

**Example 2.1**: Let $X = \{a, b\}$ then

$0_\mathcal{F} = \{(a, 0), (b, 0)\}$, $1_\mathcal{F} = \{(a, 1), (b, 1)\}$, $A_1 = \{(a, 0), (b, 0.5)\}$, $A_2 = \{(a, 0), (b, 1)\}$, $A_3 = \{(a, 1), (b, 0.2)\}$, $A_4 = \{(a, 0.1), (b, 0.5)\}$, $A_5 = \{(a, 0.4), (b, 0)\}$, etc. are fuzzy subsets of $X$ on [0,1].

Here $\tau_1 = \{0_\mathcal{F}, 1_\mathcal{F}, A_1\}$, $\tau_2 = \{0_\mathcal{F}, 1_\mathcal{F}, A_2\}$, $\tau_3 = \{0_\mathcal{F}, 1_\mathcal{F}, A_1, A_2\}$,... are fuzzy topologies on $X$.

We have trivially $\tau_\mathcal{F}(n, m, 2) = 1$ and $\tau_\mathcal{F}(n, m, 3) = m^n - 2$.

**Proposition 2.1 (M. Benoumhani [11])**: The number of fuzzy topologies on $X$, with membership values in $M$, is finite if and only if $X$ and $M$ are both finite.

**Theorem 2.1 (M. Benoumhani [11])**: The number of fuzzy topologies in $\mathcal{F}$ consisting of four open sets is given by:

$$\tau_\mathcal{F}(n, m, 4) = \left(\frac{m(m+1)}{2}\right)^n - 3m^n + 2^{n-1} + 2.$$

**Theorem 2.2 (M. Benoumhani [11])**: The number of fuzzy topologies in $\mathcal{F}$ consisting of five open sets is given by:

$$\tau_\mathcal{F}(n, m, 5) = \binom{m+2}{3}^n - 4\binom{m+1}{2}^n + (2m-1)^n + 5m^n - (m-1)^n - 2^{n+1}.$$

**Theorem 2.3 (M. Benoumhani [11])**: For $n \geq m \geq 2$, there are exactly $n(n-1)$ non-discrete topologies of maximal cardinality. Every such fuzzy topology has exactly $m^n - m^{n-2}$ open sets. That is for $n \geq m \geq 2$, we have:

(1) $\tau_\mathcal{F}(n, m, k) = 0$ for $m^n - m^{n-2} < k < m^n$, amd
(2) $\tau_\mathcal{F}(n, m, m^n - m^{n-2}) = n(n-1)$.

**Definition 2.3 [12]**: A fuzzy bitopological space is a triple $(X, \tau_1, \tau_2)$, where $\tau_1$ and $\tau_2$ are arbitrary fuzzy topologies on $X$.

**Example 2.2**: Let $X = \{a, b, c\}$ then

$0_\mathcal{F} = \{(a, 0), (b, 0), (c, 0)\}$, $1_\mathcal{F} = \{(a, 1), (b, 1), (c, 1)\}$, $A_1 = \{(a, 0), (b, 0.5), (c, 0)\}$, $A_2 = \{(a, 0), (b, 1), (c, 0)\}$, $A_3 = \{(a, 1), (b, 0.2), (c, 0.2)\}$, $A_4 = \{(a, 0.1), (b, 0.5), (c, 0)\}$, $A_5 = \{(a, 0.4), (b, 0), (c, 1)\}$, etc. are fuzzy subsets of $X$ on [0,1].

Here $\tau_1 = \{0_{\mathcal{F}}, 1_{\mathcal{F}}, A_1\}$, $\tau_2 = \{0_{\mathcal{F}}, 1_{\mathcal{F}}, A_2\}$, $\tau_3 = \{0_{\mathcal{F}}, 1_{\mathcal{F}}, A_3\}$, $\tau_4 = \{0_{\mathcal{F}}, 1_{\mathcal{F}}, A_1, A_2\}$ etc are fuzzy topologies on $X$. Then $(X, \tau_1, \tau_2)$, $(X, \tau_1, \tau_3)$, $(X, \tau_2, \tau_3)$ etc. are fuzzy bitopological spaces

### 3. Main Results:

**Result 3.1:** For any finite $n \geq 1, m \geq 2$,

(a) The number of fuzzy bitopological space in $\mathcal{F}$ consisting two open set in both the fuzzy topology is $(\tau_i, \tau_j)_{\mathcal{F}}(n, m, 2) = 1$ where $\tau_{\mathcal{F}}(n, m, 2) = 1$.

(b) The number of fuzzy bitopological space in $F$ consisting three open set in both the fuzzy topology is $(\tau_i, \tau_j)_{\mathcal{F}}(n, m, 3) = \frac{\tau_{\mathcal{F}}(n,m,3)\{\tau_{\mathcal{F}}(n,m,3)+1\}}{2} = \frac{m^{2n}-3m^n+2}{2}$, where

$\tau_{\mathcal{F}}(n, m, 3) = m^n - 2$.

**Result 3.2:** The number of fuzzy bitopological space in $\mathcal{F}$ consisting four open set in both the fuzzy topology is $(\tau_i, \tau_j)_{\mathcal{F}}(n, m, 4) = \frac{\tau_{\mathcal{F}}(n,m,4)(\tau_{\mathcal{F}}(n,m,4)+1)}{2}$, where

$\tau_{\mathcal{F}}(n, m, 4) = \left(\frac{m(m+1)}{2}\right)^n - 3m^n + 2^{n-1} + 2$.

**Example 3.2.1:** Let $X = \{a, b\}$ and $M = \{0, 0.5, 1\}$ then fuzzy subsets of $X$ with membership values in $M$ are

$0_{\mathcal{F}} = \{(a, 0), (b, 0)\}$, $1_{\mathcal{F}} = \{(a, 1), (b, 1)\}$, $A_1 = \{(a, 0), (b, 0.5)\}$, $A_2 = \{(a, 0), (b, 1)\}$,
$A_3 = \{(a, 1), (b, 0)\}$, $A_4 = \{(a, 1), (b, 0.5)\}$, $A_5 = \{(a, 0.5), (b, 0)\}$, $A_6 = \{(a, 0.5), (b, 0.5)\}$,
$A_7 = \{(a, 0.5), (b, 1)\}$

Here $n = 2, m = 3$, so

$$\tau_{\mathcal{F}}(2,3,4) = \left(\frac{3(3+1)}{2}\right)^2 - 3 \cdot 3^2 + 2^{2-1} + 2$$

$$= \left(\frac{3 \cdot 4}{2}\right)^2 - 27 + 2 + 2$$

$$= 6^2 - 23$$

$$= 36 - 23$$

$$= 13$$

These fuzzy topologies are

$\tau_1 = \{0_{\mathcal{F}}, 1_{\mathcal{F}}, A_1, A_2\}$, $\tau_2 = \{0_{\mathcal{F}}, 1_{\mathcal{F}}, A_1, A_4\}$, $\tau_3 = \{0_{\mathcal{F}}, 1_{\mathcal{F}}, A_1, A_6\}$, $\tau_4 = \{0_{\mathcal{F}}, 1_{\mathcal{F}}, A_1, A_7\}$,
$\tau_5 = \{0_{\mathcal{F}}, 1_{\mathcal{F}}, A_2, A_7\}$, $\tau_6 = \{0_{\mathcal{F}}, 1_{\mathcal{F}}, A_3, A_4\}$, $\tau_7 = \{0_{\mathcal{F}}, 1_{\mathcal{F}}, A_3, A_5\}$, $\tau_8 = \{0_{\mathcal{F}}, 1_{\mathcal{F}}, A_4, A_5\}$,
$\tau_9 = \{0_{\mathcal{F}}, 1_{\mathcal{F}}, A_4, A_6\}$, $\tau_{10} = \{0_{\mathcal{F}}, 1_{\mathcal{F}}, A_5, A_6\}$, $\tau_{11} = \{0_{\mathcal{F}}, 1_{\mathcal{F}}, A_5, A_7\}$, $\tau_{12} = \{0_{\mathcal{F}}, 1_{\mathcal{F}}, A_6, A_7\}$,
$\tau_{13} = \{0_{\mathcal{F}}, 1_{\mathcal{F}}, A_2, A_3\}$.

So, the number of fuzzy bitopological space is

$$(\tau_i, \tau_j)_{\mathcal{F}}(2,3,4) = \frac{T_{\mathcal{F}}(2,3,4)(T_{\mathcal{F}}(2,3,4) + 1)}{2}$$

$$= \frac{13(13+1)}{2}$$

$$= \frac{13 \cdot 14}{2}$$

$$= 13 \times 7$$

$$= 91.$$

**Result 3.3:** The number of fuzzy bitopologies in $\mathcal{F}$ consisting of five open sets in both the fuzzy topology is $(\tau_i, \tau_j)_{\mathcal{F}}(n, m, 5) = \frac{T_{\mathcal{F}}(n,m,5)(T_{\mathcal{F}}(n,m,5)+1)}{2}$, where

$$T_{\mathcal{F}}(n, m, 5) = \binom{m+2}{3}^n - 4\binom{m+1}{2}^n + (2m-1)^n + 5m^n - (m-1)^n - 2^{n+1}.$$

**Example 3.3.1:** Let $X = \{a, b\}$ and $M = \{0, 0.2, 1\}$ then fuzzy subsets of $X$ with membership values in $M$ are

$0_{\mathcal{F}} = \{(a, 0), (b, 0)\}$, $1_{\mathcal{F}} = \{(a, 1), (b, 1)\}$, $A_1 = \{(a, 0), (b, 0.2)\}$, $A_2 = \{(a, 0), (b, 1)\}$,
$A_3 = \{(a, 0.2), (b, 0)\}$, $A_4 = \{(a, 0.2), (b, 1)\}$, $A_5 = \{(a, 0.2), (b, 0.2)\}$,
$A_6 = \{(a, 1), (b, 0.2)\}$, $A_7 = \{(a, 1), (b, 0)\}$.

Here $n = 2, m = 3$ so

$$T_{\mathcal{F}}(2,3,5) = \binom{3+2}{3}^2 - 4\binom{3+1}{2}^2 + (2 \cdot 3 - 1)^2 + 5 \cdot 3^2 - (3-1)^2 - 2^{2+1}$$

$$= \binom{5}{2}^2 - 4 \cdot \binom{4}{2}^2 + 25 + 45 - 4 - 8$$

$$= \left(\frac{5!}{2!3!}\right)^2 - 4 \cdot \left(\frac{4!}{2!2!}\right)^2 + 70 - 12$$

$$= 100 - 144 + 58$$

$$= 14$$

These fuzzy topologies are

$\tau_1 = \{0_{\mathcal{F}}, 1_{\mathcal{F}}, A_1, A_2, A_4\}$, $\tau_2 = \{0_{\mathcal{F}}, 1_{\mathcal{F}}, A_1, A_2, A_7\}$, $\tau_3 = \{0_{\mathcal{F}}, 1_{\mathcal{F}}, A_1, A_3, A_4\}$,
$\tau_4 = \{0_{\mathcal{F}}, 1_{\mathcal{F}}, A_1, A_3, A_7\}$, $\tau_5 = \{0_{\mathcal{F}}, 1_{\mathcal{F}}, A_1, A_4, A_6\}$, $\tau_6 = \{0_{\mathcal{F}}, 1_{\mathcal{F}}, A_1, A_5, A_6\}$,
$\tau_7 = \{0_{\mathcal{F}}, 1_{\mathcal{F}}, A_1, A_5, A_7\}$, $\tau_8 = \{0_{\mathcal{F}}, 1_{\mathcal{F}}, A_1, A_6, A_7\}$, $\tau_9 = \{0_{\mathcal{F}}, 1_{\mathcal{F}}, A_2, A_5, A_7\}$,
$\tau_{10} = \{0_{\mathcal{F}}, 1_{\mathcal{F}}, A_3, A_4, A_5\}$, $\tau_{11} = \{0_{\mathcal{F}}, 1_{\mathcal{F}}, A_3, A_5, A_7\}$, $\tau_{12} = \{0_{\mathcal{F}}, 1_{\mathcal{F}}, A_4, A_5, A_6\}$,
$\tau_{13} = \{0_{\mathcal{F}}, 1_{\mathcal{F}}, A_4, A_6, A_7\}$, $\tau_{14} = \{0_{\mathcal{F}}, 1_{\mathcal{F}}, A_5, A_6, A_7\}$.

So, the number of fuzzy bitopological space is

$$(\tau_i, \tau_j)_{\mathcal{F}}(2,3,5) = \frac{\tau_{\mathcal{F}}(2,3,5)(\tau_{\mathcal{F}}(2,3,5)+1)}{2}$$

$$= \frac{14(14+1)}{2}$$

$$= \frac{14 \cdot 15}{2}$$

$$= 7 \times 15$$

$$= 105.$$

**Result 3.4:** The number of fuzzy bitopological space having $k$ open set in both fuzzy topology

(a) $(\tau_i, \tau_j)_{\mathcal{F}}(n, m, k) = 0$ for $m^n - m^{n-2} < k < m^n$ and

(b) $(\tau_i, \tau_j)_{\mathcal{F}}(n, m, k) = \frac{\tau_{\mathcal{F}}(n,m,k)\{\tau_{\mathcal{F}}(n,m,k)+1\}}{2} = \frac{n^4 - 2n^3 + 2n^2 - n}{2}$ where $k = m^n - m^{n-2}$ and $\tau_{\mathcal{F}}(n, m, k) = n(n-1)$.

## 4. Conclusion:

If we can determine number of fuzzy topologies having $k$ open set $\tau_{\mathcal{F}}(n, m, k)$ on a finite set then we can also determine number of fuzzy bitopological space $(\tau_i, \tau_j)_{\mathcal{F}}(n, m, k)$.